%%%% The simple exclusion process with multiple shocks %%%%%%%%%
%%%% by Ferrari, Fontes and Vares %%%%%%%%
%%%% 02\05/99 r,e %%%%%%%%
%%%% final
%%%%%%%%%%%%%%%%%%%
\magnification=1200
\input amstex
\input amsppt.sty
%\nologo
\openup 6pt
\NoBlackBoxes

%\hsize=4.9 true in
%\vsize=9.0 true in
\vsize=22.5truecm
%\hoffset=-.250 true in
%\voffset=-.150 true in

%%%%%%%%%%%%%%%%% EQUAZIONI CON NOMI SIMBOLICI
%%%
%%% Per assegnare un nome simbolico ad una equazione basta
%%% scrivere \Eq(...) o, in \eqalignno, \eq(...) o,
%%% nelle appendici, \Eqa(...) o \eqa(...);
%%% dentro le parentesi e al posto di ... si puo' scrivere qualsiasi commento;
%%% per avere i nomi simbolici segnati a sinistra delle formule si deve
%%% dichiarare il documento come bozza, iniziando il testo con
%%% \BOZZA. Sinonimi: \Eq,\EQ,\EQS; \eq,\eqs; \Eqa,\Eqas;\eqa,\eqas.
%%% All' inizio di ogni paragrafo si devono definire il
%%% numero del paragrafo e della prima formula dichiarando
%%% \numsec=... \numfor=...  (brevetto Eckmannn).
%%% Si possono citare formule seguenti; le corrispondenze fra nomi
%%% simbolici e numeri effettivi sono memorizzate nel file \jobname.aux, che
%%% viene letto all'inizio, se gia' presente. E' possibile citare anche
%%% formule che appaiono in altri file, purche' sia presente il
%%% corrispondente file .aux; basta includere all'inizio l'istruzione
%%%           \include{nomefile}
%%%
%%%%%%%%%%%%%%%%%%%%%%%%%%%%%%%%%%%%%%%%%%%%%%%%%%%%%%%%%%%%%%%

\global\newcount\numsec
\global\newcount\numfor
\global\newcount\numtheo
\global\advance\numtheo by 1

\def\senondefinito#1{\expandafter\ifx\csname#1\endcsname\relax}

\def\SIA #1,#2,#3 {\senondefinito{#1#2}
\expandafter\xdef\csname #1#2\endcsname{#3}\else
\write16{???? ma #1,#2 e' gia' stato definito !!!!} \fi}

\def\etichetta(#1){(\veroparagrafo.\veraformula)%
\SIA e,#1,(\veroparagrafo.\veraformula)
\global\advance\numfor by 1
\write15{\string\FU (#1){\equ(#1)}}%
\write16{ EQ #1 ==> \equ(#1) }}%

\def\letichetta(#1){\veroparagrafo.\verotheo%  thms etc keep own numeration
\SIA e,#1,{\veroparagrafo.\verotheo}
\global\advance\numtheo by 1
\write15{\string\FU (#1){\equ(#1)}}%
\write16{ Sta \equ(#1) == #1 }}%

\def\tetichetta(#1){\veroparagrafo.\veraformula%%%%copy four lines
\SIA e,#1,{(\veroparagrafo.\veraformula)}
\global\advance\numfor by 1%
\write15{\string\FU (#1){\equ(#1)}}%
\write16{ tag #1 ==> \equ(#1)}}%

\def\FU(#1)#2{\SIA fu,#1,#2 }

\def\etichettaa(#1){(A\veroparagrafo.\veraformula)%
\SIA e,#1,(A\veroparagrafo.\veraformula) %
\global\advance\numfor by 1%
\write15{\string\FU (#1){\equ(#1)}}%
\write16{ EQ #1 ==> \equ(#1) }}

\def\BOZZA{
\def\alato(##1){%
 {\rlap{\kern-\hsize\kern-1.4truecm{$\scriptstyle##1$}}}}%
\def\aolado(##1){%
{%\vtop to \profonditastruttura
{%\baselineskip
%\profonditastruttura\vss
\rlap{\kern-1.4truecm{$\scriptstyle##1$}}}}}%
}

\def\alato(#1){}
\def\aolado(#1){}

\def\veroparagrafo{\number\numsec}
\def\veraformula{\number\numfor}
\def\verotheo{\number\numtheo}

\def\Eq(#1){\eqno{\etichetta(#1)\alato(#1)}}
\def\eq(#1){\etichetta(#1)\alato(#1)}
\def\Leq(#1){\leqno{\aolado(#1)\tetichetta(#1)}$\!\!$}%%%%%this line for \leqno
\def\teq(#1){\tag{\aolado(#1)\tetichetta(#1)\alato(#1)}$\!$}%%%%%this line for \tag
\def\Eqa(#1){\eqno{\etichettaa(#1)\alato(#1)}}
\def\eqa(#1){\etichettaa(#1)\alato(#1)}
\def\eqv(#1){\senondefinito{fu#1}$\clubsuit$#1
\write16{#1 non e' (ancora) definito}%
\else\csname fu#1\endcsname\fi}
\def\equ(#1){\senondefinito{e#1}\eqv(#1)\else\csname e#1\endcsname\fi}

%next six lines by paf (no responsibilities taken)
\def\Lemma(#1){\aolado(#1)Lemma \letichetta(#1)$\!\!$}%
\def\Theorem(#1){{\aolado(#1)Theorem \letichetta(#1)}$\!\!$}%\hskip-1.6truemm}%
\def\Proposition(#1){\aolado(#1){Proposition \letichetta(#1)}$\!\!$\hskip-1.6truemm}%
\def\Corollary(#1){{\aolado(#1)Corollary \letichetta(#1)}$\!\!$\hskip-1.6truemm}%
\def\Remark(#1){{\noindent\aolado(#1){\bf Remark \letichetta(#1)$\!$.}}}%
\def\Definition(#1){{\noindent\aolado(#1){\bf Definition \letichetta(#1)$\!\!$\hskip-1.6truemm}}}%
\def\Example(#1){\aolado(#1) Example \letichetta(#1)$\!\!$\hskip-1.6truemm}
\let\ppclaim=\plainproclaim

\def\include#1{
\openin13=#1.aux \ifeof13 \relax \else
\input #1.aux \closein13 \fi}

\openin14=\jobname.aux \ifeof14 \relax \else
\input \jobname.aux \closein14 \fi
\openout15=\jobname.aux

%\BOZZA
%%%%%%%%%%%%%%%%%   END  NUME2.TEX   %%%%%%%%%%%%%%%%%%%%%%%%%%%%%%%%%%%%%
%%%%%%%%%%%%%%%%%%%%%%%%%%%%%%%%%%%%%%%%%%%%%%%%%%%%%%%%%%%%%%%%%%%%%%%%%%

\font\bigbf=cmbx10 scaled \magstep1

\def\Z{\Bbb Z}
\def\r{\Bbb R}
\def\ro{\rho}

\def\be{\beta}

\def\ve{\varepsilon}
\def\de{\delta}

\def\Om{\Omega}
\def\om{\omega}
\def\a{\alpha}
\def\aa{\Cal A}

\def\x{\Cal X}
\def\y{\Cal Y}
\def\n{\Cal N}

\def\t{\theta}

\def\ov{\overline}
\def\un{\underline}
\def\po{\partial}

\def\lv{\left\vert }
\def\rv{\right\vert }

\def\ga{\gamma}
\def\si{\sigma}

\def\La{\Lambda}
\def\={&=}

\def\one{{\bold 1}}
\def\square{\ifmmode\sqr\else{$\sqr$}\fi} 
\def\sqr{\vcenter{ 
         \hrule height.1mm 
         \hbox{\vrule width.1mm height2.2mm\kern2.18mm\vrule width.1mm} 
         \hrule height.1mm}}                  % This is a slimmer sqr. 
\def\ie{{\it i.e.}\ }

\centerline{\bigbf THE ASYMMETRIC SIMPLE EXCLUSION PROCESS}
\centerline{{\bigbf WITH MULTIPLE SHOCKS}}

\vglue .5in

\centerline{Pablo A. Ferrari${}^{(1)}$, L. Renato G. Fontes${}^{(1)}$,
M. Eul\'alia Vares${}^{(2)}$}

\footnote""{Partially Supported by CNPq and FAPESP and FINEP (PRONEX).\newline
${}^{(1)}$ IME-USP. Caixa Postal 66.281, S\~ ao Paulo 05315-970 SP
Brasil\newline 
${}^{(2)}$ IMPA. Estrada D. Castorina 110, J. Bot\^anico 22460-320 RJ Brasil}

%%%%%%%%%%%%%%%%%%%%%%%%%%%%%%%%%%%%%%%%%%%%%%%%%%%%%%%%%%%%%%%%%%%%%%%%%

\vglue .3in

\baselineskip 10pt \flushpar {\bf Abstract:} We consider the one
dimensional totally asymmetric simple exclusion process with initial
product distribution with densities $0 \leq \rho_0 < \rho_1 <...<
\rho_n \leq 1$ in $(-\infty,c_1\ve^{-1})$,
$[c_1\ve^{-1},c_2\varepsilon^{-1}),\dots,[c_n \ve^{-1} , + \infty)$,
respectively. The initial distribution has shocks (discontinuities) at
$\varepsilon^{-1}c_k$, $k=1,\dots,n$ and we assume that in the
corresponding macroscopic Burgers equation the $n$ shocks meet in
$r^*$ at time $t^*$. The microscopic position of the shocks is
represented by second class particles whose distribution in the scale
$\varepsilon^{-1/2}$ is shown to converge to a function of $n$
independent Gaussian random variables representing the fluctuations of
these particles ``just before the meeting''. We show that the density
field at time $\ve^{-1}t^*$, in the scale $\ve^{-1/2}$ and as seen
from $\ve^{-1}r^*$ converges weakly to a random measure with piecewise
constant density as $\ve \to 0$; the points of discontinuity depend on
these limiting Gaussian variables. As a corollary we show that, as
$\varepsilon\to 0$, the distribution of the process at site
$\varepsilon^{-1}r^*+\ve^{-1/2}a$ at time $\varepsilon^{-1}t^*$ tends
to a non trivial convex combination of the product measures with
densities $\rho_k$, the weights of the combination being explicitly
computable.

\vskip 3mm

shock fluctuations

\vskip 3mm 

\vskip  3mm

\vfill\eject
%%%%%%%%%%%%%%%%%%%%%%%%%%%%%%%%%%%%%%%%%%%%%%%%%%%%%%%%%%%%%%%%%%%%%%%%%

%\vglue .5in

\flushpar
{\bf \S 1.\quad Introduction and results.}
\numsec=1\numfor=1\numtheo=1
\bigskip

It is well known that the hydrodynamical behavior of the one dimensional
asymmetric exclusion process is described by the inviscid Burgers equation
$$
\po_t\ro+\ga\po_r(\ro(1-\ro))=0\teq(1.1)
$$
where $\ga$ is the mean of the jump distribution. Since equation (1.1) develops
discontinuities, one has to be careful about the precise statement, but loosely
speaking, if $(r,t)$ is a continuity point of $\ro(r,t)$, for a given initial
measurable profile $\ro_0(\,\,\cdot\,\,)$, then at the macroscopic point
$(r,t)$ the system is distributed according to the measure
$\nu_{\ro(r,t)}$, where $\nu_\ro$ is the product Bernoulli measure on
$\{0,1\}^{\Z}$ with $\nu_\ro(\eta(x)=1)=\ro$ for all $x$. 
This is known as local equilibrium.
As it is known,
the exact statement involves a space/time change, under Euler scale, and for
all these developments we refer to Andjel and Vares (1987), Rezakhanlou
(1990), Landim (1992).

The problem with which we are concerned here involves the description
of the (microscopic) behavior of the system at certain discontinuity
points $(r,t)$ (or shock fronts) of the solution of equation
(1.1). For example, if $\ga>0$ and the initial profile is a step
function $\ro_0(r):=\a\bold 1\{r<0\}+ \be\bold 1\{r\ge0\}$, with
$0\le\a<\be\le 1$, the entropy solution of equation (1.1) is
$\ro(r,t)=\a\bold 1\{r<vt\}+\be\bold 1\{r\ge vt\}$, where
$v:=\ga(1-\a-\be)$ is the velocity of the shock front and
$\one\{\cdot\}$ is the indicator function of the set $\{\cdot\}$.
This description is valid only for continuity points.  The
investigation of what happens to the system if looked from this shock
front was first studied by Wick (1985) for a different model, and for
the asymmetric simple exclusion in the particular situations $\a=0$
and $\a+\be=1$ by De Masi et al (1988) and Andjel, Bramson and Liggett
(1988), respectively. They all proved that at the shock front one sees
a fair mixture of $\nu_\a$ and $\nu_\be$. This result has then been
extended so as to cover all cases $0\le\a<\be\le1$ by Ferrari and
Fontes (1994), from now on referred as [FF].  [FF] worked with the
nearest neighbor asymmetric exclusion process, whose generator is the
closure of
$$
Lf(\eta):=\sum_{x\in\Z}\sum_{y=x\pm 1} p(x,y)\eta(x)(1-\eta(y))(f(\eta^{x,y})-
f(\eta))\teq(1.2)
$$
for $f$ a cylinder function in $\{0,1\}^{\Z}$, with
$$
\eta^{x,y}(z):=\cases
\eta(z) & \text{ if } z\ne x,y\\
\eta(y) & \text{ if } z=x \\
\eta(x) & \text{ if } z=y, \endcases
$$
where $p(x,x+1):=p$, $p(x,x-1):=q:=1-p$, with $1/2<p\le1$. This process was first
studied by Spitzer (1970).
Calling $\mu_{\a,\be}$ the product measure on
$\{0,1\}^{\Z}$ with site marginals
$$
\mu_{\a,\be}(\eta(x)=1):= \cases
\a &\text{ if } x<0\\
\be & \text{ if } x\ge0, \endcases\teq(1.3)
$$ 
denoting $S_t$ as the semigroup corresponding to the above
generator and $\t_x$ as the space shift $\mu \theta_x(f):= \int
f(\theta_x \eta) \mu(d\eta)$, with $\t_x\eta(z):=\eta(x+z)$, [FF]
proved that
$$
\mu_{\a,\be} S_t\t_{[vt]}\overset{\om^*}\to{\longrightarrow} \frac12
(\nu_\a+\nu_\be)\teq(1.4)
$$ 
as $t\to+\infty$, and where $[x]$ denotes the integer part of $x$. 
This corresponds to the exact statement made
above, under Euler scale, and for the macroscopic point $(r,t)$ in the
front line, \ie $r=vt$.  In fact, in the above mentioned references,
more detailed analysis is performed, by looking at the microscopic
structure of the shock represented by a second class particle. A
second class particle jumps to empty sites with the same rates as the
other particles, but interchanges positions with the regular particles
at the rate holes do. A formal definition using coupling is given in
the next section. Calling $X_t$ the position of a second class
particle added at the origin, the process as seen from the second
class particle $\theta_{X_t}\eta_t$ has distribution 
asymptotically
product to the right and left of the origin with densities $\alpha$
and $\beta$ respectively, uniformly in time. The velocity of the second
class particle is the same as the velocity of the shock in the Burgers
equation: $E_{\mu_{\alpha,\beta}} X_t = \gamma(1-\alpha-\beta)t$.  [FF] proved
that the fluctuations of the position of the particle are Gaussian:
calling $\tilde X_t := X_t - vt$, 
$$
\tilde X_t/\sqrt t\ \ \overset{\Cal D}\to{\underset{t \to +
\infty}\to{\longrightarrow}}\ \  \n_{\alpha,\beta}
\teq(p1)
$$
where $\n_{\alpha,\beta}$ is a centered Gaussian random
variable with variance $\gamma(\beta-\alpha)^{-1}
(\alpha(1-\alpha)+\beta(1-\beta))$.  With this result in hand [FF]
proved that, if $ -\infty < a < +\infty$, the distribution of the
process at time $t$ at the point $vt+a\sqrt t$ converges to a mixture
of $\nu_\a$ and $\nu_\be$; more precisely, for real $a$,
$$ 
\mu_{\a,\be} S_t\t_{[vt+a\sqrt t]}\ \
\overset{\om^*}\to{\underset{t \to +\infty}\to{\longrightarrow}}\ \ \nu_\a P(\n_{\alpha,\beta} >a)
+\nu_\be P(\n_{\alpha,\beta} \le a)\teq(p2)
$$ 
which in particular yields $\frac12 (\nu_\alpha + \nu_\beta)$  for  $a=0$
and interpolates between $\nu_\alpha$ and $\nu_\beta$, as $a$ varies from
$-\infty$ to $+\infty$. 

Our goal is the consideration of two or more shock fronts and the
description of the microscopic behavior of the system at the
(macroscopic) time and position of their meeting. To avoid unnecessary
technical difficulties, we look the totally asymmetric case: we assume
$$
p=1,\,\, q=0\,\,\,\,\text{ which implies }\,\,\,\,
\gamma=1\,\,\,\,\text{ in equation }\equ(1.1). 
$$

The extension to $\frac12<p<1$ will be briefly
discussed at the end.  

We consider points $c_k$, densities $\rho_k$ and the existence of a
space-time point $ (r^*,t^*$) such that
$$
\align
&-\infty=c_0<c_1<\dots<c_n<c_{n+1}=\infty, \teq(p10a)\cr
&0\leq\rho_0<\dots<\rho_{n}\leq1, \teq(p10b)\cr
&r^* = c_k + (1-\rho_{k-1}-\rho_{k}) t^*,\quad k=1,\dots,n.\teq(p10)\cr
\endalign
$$
With this assumption the entropy solution to equation \equ(1.1) with initial
data
$$
\ro(r)\,:=\, \sum_{k=0}^{n} \ro_k\,\bold 1\{c_{k}\le r <c_{k+1})\},\teq(1.5)
$$ 
has the property that all $n$ shocks meet at $r^*$ at time $t^*$.
More precisely, the entropy solution is given by
$$
\ro(r,t)\,=\, \sum_{k=0}^n \ro_k\,
\bold 1\{c_{k}(t)\le r <c_{k+1}(t))\}\teq(1.6b)
$$
where 
$$
c_k(t) = 
\cases c_k + (1-\rho_{k-1}-\rho_{k}) t &\text{ for } t<t^* \cr
r^* + (1-\rho_0-\rho_n)(t-t^*) &\text{ for } t\ge t^* \cr
\endcases \teq(1.6p)
$$ 
Notice that after $t^*$ only the extreme densities $\rho_0$ and
$\rho_n$ are seen. 

In our discussion we shall need to consider entropy solutions starting
with more general increasing step profiles, {\ie not
necessarily with all shocks meeting at $(t^*,r^*)$. For this let
$(\rho_k)$ satisfy \equ(p10b), $-\infty =b_0< b_1 < ...<b_n <b_{n+1}=
+\infty$ and consider the entropy solution to the Burgers equation
with initial data
$$ 
\lambda(r)\,:=\, \sum_{k=0}^{n} \ro_k\,\bold 1\{b_{k}\le r
<b_{k+1})\},\teq(p15)
$$ 
which is given by
$$
\align  
\lambda(r,t)\,&=\, \sum_{k=0}^n \ro_k\,
\bold 1\{b_{k}(t)\le r <b_{k+1}(t))\}\teq(p16)\cr
\displaystyle {d b_k(t) \over dt} \,&=\, 1- \lambda^+(b_k(t),t) -
\lambda^-(b_k(t),t)  \teq(p12)\cr
 b_k(0) \,&=\, b_k,\;\; i=1,\dots,n. 
\endalign
$$ 
where $\lambda^\pm(r,t)$ are the right and left limits of
$\lambda(r,t)$: $\lambda^\pm(r,t) = \lim_{r'\to0}\lambda(r\pm r',t)$,
with $r'>0$.  In words, $b_1(t) \leq \dots\leq b_n(t)$ represent the
shock fronts, which move initially as $b_j(t)=b_j + t(1-\rho_{j-1} -
\rho_j)$ for $1 \leq j \leq n$, until two or more of them meet. When
this happens the involved shock fronts coalesce, with the
disappearance of all intermediate densities and the front keeps moving
with a new velocity given by one minus the sum of the two densities, to
its left and to its right, \ie the two densities which form the
shock.  This simple description, mainly due to the fact that we are in
the one dimensional situation and the initial profile is an increasing
step function, allows to define the following map $\psi:\r^n\to\r^n$,
$\psi= (\psi_k)$, which will be fundamental in the determination of
the distribution of the process at the macroscopic meeting point of
the shocks:

\noindent{\bf Definition of $\psi$. } 
Given a vector $\underline x =(x_1,\dots,x_n)$ in $ \Bbb R^n$ let us take a
time $t(\underline x)$ large enough so that defining 
$$ 
 b_k(\underline x)\, :=\, x_k - t(\underline
x)\,(1-\rho_{k-1}-\rho_k)\teq(p222)
$$
then 
$$
b_1(\underline x)<\dots<b_n(\underline x).\teq(p223)
$$ 
That is, if we consider a family of one shock solutions
$\lambda^k(r,t) ,\; 1 \leq k \leq n$ starting with   
$$
\lambda^k(r):= \rho_{k-1} \text{\bf 1}\{r < b_k(\underline x)\} +
\rho_k \text{ \bf 1} \{r \geq b_k( \underline x)\} 
$$ 
then $x_k$ is the position of the shock at time $t(\un x)$, for the
$k$-th equation. What we denote by $\psi(\underline x)$ is then the
position of the shock fronts at time $t(\un x)$ in equations
\equ(p16)-\equ(p12) starting with \equ(p15) for $(b_k
=b_k(\underline x))$. It is easy to see that the definition is
well posed, \ie it does not depend on the value of $t(\underline
x)$ provided \equ(p222) and \equ(p223) hold.

The coordinates of the vector $\psi(\un x)= (\psi_k(\underline x))$
are convex combinations of some $x_j$. In the particular case of $n=2$
this becomes
$$
\psi_k(x_1,x_2) = \cases x_k &\text{if } x_1\le x_2\cr
\displaystyle{x_1 {\rho_1-\rho_0\over \rho_2-\rho_0} + x_2
{\rho_2-\rho_1\over \rho_2-\rho_0}} &\text{if } x_1>x_2\cr
\endcases\teq(p28)
$$ 
for $k=1,2$. 

We consider a family of product measures $\mu^\ve$ on $\{0,1\}^\Z$
with profile $\rho$: the marginal one-site distribution is given by
$$
\mu^\ve(\eta(x)=1) = \rho(\ve x),
$$ 
where $\rho(\cdot)$ is given in \equ(1.5) and corresponds to the case
when all shocks meet at the same point $r^*$ at time $t^*$. Considering
$\mu^\ve$ as the initial measure, 
our first goal is to look at the asymptotic distribution of the
process at time $\tau_\ve$ around site 
$[x_\ve]$ in the scale $\ve^{-1/2}$ as $\ve \to 0$,  where 
$$
\tau_\ve:=t^*\ve^{-1},\;\;\;\;\;\; x_\ve :=r^*\ve^{-1}.\teq(tre) 
$$

Let $\y_0:=-\infty$, $\y_{n+1}:=\infty$ and 
$$
(\y_1,\dots,\y_n):=
\psi(\x_1,\dots,\x_n)\teq(p9)
$$
where $(\x_1,\dots,\x_n)$  are independent
centered Gaussian random variables with variances $(D_1,\dots,D_n)$
given by
$$
D_k := {\rho_{k-1}(1-\rho_{k-1})+\rho_{k}(1-\rho_{k})\over
\rho_{k}-\rho_{k-1}}\,t^*\teq(p99)
$$

\bigskip

\ppclaim \Theorem(t1.1). Let $\mu^\ve$ $(0<\ve\le 1)$ be the
product measure on $\{0,1\}^{\Z}$, associated to the profile
$\ro(\,\cdot\,)$, given by \equ(1.5). Let $S_t$ denote the
semigroup associated to the totally asymmetric $n.n.$ simple exclusion
process, corresponding to the generator given by \equ(1.2) with
$p=1$. Then, for $a\in \r$,
$$
\mu^\ve S_{t^*\ve^{-1}}\theta_{[r^*\ve^{-
1}+a\ve^{-1/2}]}\ \overset{\om^*}\to{\underset{\ve\to0}\to{\longrightarrow}}\
\sum_{k=0}^{n}    \nu_{\rho_k}P(\y_{k} \le a < \y_{k+1}) \teq(p5)
$$ 
where $\y_k$ are defined in \equ(p9). 

\bigskip 

We shall see that, as in \equ (p1), $(\Cal Y_1,..., \Cal Y_n)$ represents the limiting
fluctuations of the shocks at time $t^*\ve^{-1}$, around 
$r^*\ve^{-1}$. In the case of only one shock $(n=1)$ this agrees with
\equ(p2). The crucial difference is that if, say, the $k$-th and
$(k+1)$-th microscopic shocks meet before $t^*\ve^{-1}$,
they coalesce and change the velocity and the intermediate zone of
density $\rho_k$ disappears. This explains the weight of $\nu_{\rho_k}$in \equ(p5):
it is the same as the probability that (a)
the $k$-th and $(k+1)$-th shocks have not collided yet and (b) these
shocks are to the left and right of the point we are looking at,
respectively. 
%Since they are convex combination of normal random
%variables of variance of order $\ve^{-1}$, they will be far away of
%the point we are looking at.  
The coalescing dynamics of the
microscopic shocks in the scale $\ve^{-1/2}$ is as in the Burgers
equation. Its relation with the $n$ one-shock dynamics ---represented
by $(\x_k)$--- is given by the function $\psi$.

We turn now to the profile seen from the meeting point and time of the
shock fronts, scaled by $\varepsilon^{-1/2}$. For $\ve>0$ and a local
function $f$, consider the (random) measures on the real line given by
$$
\La_\ve(da):=\ve^{1/2}\sum_{x\in\Z}f(\theta_{x+[x_\ve]}\eta_{\tau_\ve})\,
\delta_{\ve^{1/2}x}(da),\teq(5.1)
$$
where $x_\ve$ and $\tau_\ve$
are given by \equ(tre), $\delta_{\ve^{1/2}x}$ is the Dirac delta measure at $\ve^{1/2}x$, and
$\eta_t$ is the configuration at time t for the initial measure $\mu^\ve$ . 

Let $(\y_1,\dots,\y_n)$ be as in \equ(p9) and consider the random
measure 
$$
\eqalign{\La(da)&:=
\sum_{k=0}^n \nu_{\rho_k}(f)\,\,\one\{\y_k\leq a<\y_{k+1}\}\,da\cr}\teq(5.2)
$$
\bigskip

\ppclaim \Theorem (t1.2). $\La_\ve$ converges in law to  $\La$ as $\ve\to0$,
for the usual weak topology on the space of measures.

\bigskip

The analysis is based on the well known strategy of identifying
microscopically the shocks with second class particles, defined
through the so called basic coupling of different versions of the
process starting with measures with different uniform densities. From
the dependence of their locations on the initial condition, which is a
one-shock fact as in [FF], we can ascertain their
distribution around the macroscopic meeting place $[x_\ve]$ at the
macroscopic meeting time $\tau_\ve$; in the scale
$\ve^{-1/2}$ this is given by $(\x_1,\dots,\x_n)$. If we look at the positions of the shocks at
time $\tau_\ve-\alpha\ve^{-1/2}$, the law of large numbers apply. This is also a one-shock fact
(Ferrari (1992)). If $\alpha$ is big enough these positions are
ordered, with very large probability, and they are represented by 
the variables $b_k(\x_1,\dots,\x_n)$
as in \equ(p222), taking $t(\underline x)= \alpha \ve^{-1/2}$. From 
then on, using the invariance of the product
measures involved, we can follow their trajectory up to the meeting
time through successive applications of the one-shock law of large
numbers. The final positions of the shocks in the scale $\ve^{-1/2}$
are given by the variables $\y_k + [x_\ve]$. 

Theorem \equ(t1.2) follows from the weak convergence of a suitable function of the second
class particles to the variables $\y_k$, and the asymptotic properties
of the measure as seen from the second class particles. We show
Theorem \equ(t1.1) as a corollary of Theorem \equ(t1.2), by using the
attractiveness of the process. The proof avoids proving firstly the 
translation invariance of the weak limits of Theorem \equ(t1.1), as in
the proof given by [FF] for the one-shock case.

{\bf Acknowledgment:} The authors are indebted to an anonymous referee for
bringing in the discussion of Theorem \equ(t1.2) and for noticing that
our first proof of Theorem \equ(t1.1) was indeed yielding this profile
description.

%%%%%%%%%%%%%%%%%%%%%%%%%%%%%%%%%%%%%%%%%%%%%%%%%%%%%%%%%%%%%%%%%%%%%%%%%

\vglue .3in

\flushpar {\bf \S 2.\quad Second class particles and more.}
\nopagebreak \numsec=2\numfor=1\numtheo=1 

\bigskip 

Ferrari, Kipnis and Saada (1991) and Ferrari (1992) have shown that in
the case that $\mu_{\a,\be}$ is the initial measure, the shock front
is well described by $X_t$, the position of a ``second class
particle'' initially located at the origin. This fact together with
the validity of a central limit theorem for $X_t$ as $t\to+\infty$,
proven by [FF], are the essentials for the proof that in the totally
asymmetric case $\mu_{\a,\be}S_t\theta_{[(1-\a-\be)t]}$ tends to
$\frac12(\nu_\a+\nu_\be)$. For all this, as well as in the present
work, the main tool is coupling. To realize a coupling of several
evolutions of the simple exclusion process, corresponding to several
initial configurations, is particularly simple through the graphical
construction: to each pair of sites $(x,x+1)$, let us associate a
Poisson process with rate~$1$, and at each of its occurrences we put
an arrow $x\to x+1$. Construct all such Poisson processes as
independent in some space $(\Om,\aa,P)$. Given any realization of
arrows $\om$ and a initial configuration $\eta$, we may realize an
evolution $\eta_t$ corresponding to $L$ imposing that whenever an
arrow $x\to x+1$ appears, if there is a particle at $x$ and no
particle at $x+1$, then this particle moves to $x+1$; otherwise,
nothing happens. (In the non totally asymmetric cases, we have arrows
$(x,x+1)$ with rate $p$ and $(x,x-1)$ with rate $q$, all Poisson
processes being independent.) Since the probability of two
simultaneous arrows is zero this construction makes sense and defines
the process of interest. (One could also state the coupling via a
suitable generator, cf. Liggett (1976) and Chapter VIII of Liggett
(1985).) From this coupling, the attractiveness property of the
dynamics is immediate: if $\eta^0$ and $\eta^1$ are two initial
configurations such that $\eta^0\le \eta^1$ (\ie
$\eta^0(x)\le\eta^1(x)$ $\forall x)$ and we write $\eta_t^0$,
$\eta_t^1$ for their evolutions using the same arrows, then
$\eta_t^0\le\eta_t^1$ for each $t$.

Let us then couple in this way realizations of the exclusion process
with random initial configurations, distributed according to
$\nu_{\ro_k}$, $k=0,\dots,n$. In fact we have a small perturbation
since we will add particles at $[\ve^{-1}c_k]$. For example, let $U=(U_x)_{x\in\Z}$ be a
family of i.i.d. random variables uniformly distributed on $[0,1]$,
taken as independent of all the Poisson processes of arrows, where, if
needed, we enlarge the basic probability space $(\Om,\aa,P)$. Then,
for $(c_k)$ and $(\rho_k)$ as in Section 1 define
$$
\aligned
\sigma^k([\ve^{-1}c_j]) &= \one\{j\le k\},\quad j=1,\dots,n,\cr
\sigma^k(x) &= \one\{U_x<\rho_k\},\;\;\;
x\in\Z\setminus\{[\ve^{-1}c_1],\dots,[\ve^{-1}c_n]\} . \endaligned
$$

Using the same graphical construction above described we may consider
the simultaneous evolution of all these configurations, which we
denote by $\sigma^k_t$, $k=0,\dots,n$ on the space $(\Om,\aa,P)$. The
marginal distribution of $\sigma^k_t$ is the simple exclusion process
under the invariant distribution $\nu_{\rho_k}$.

Consider the configurations on $\{0,1\}^{\Z}$ given by
$\xi^0_t=\si^0_t$ and for  $k=1,\dots,n$, 
$$
\xi^k_t(x) = \sigma^{k}_t(x)-\sigma^{k-1}_t(x)\teq(2.2a)
$$

It is easy to see that when considering the joint motion of
$(\xi^1_t,\dots,\xi^n_t)$, then for $j<k$ the $\xi^{j}$ particles have
priority over the $\xi^k$ particles: if there is a $\xi^{j}$
particle at site $x$, a $\xi^k$ particle at $x+1$ and an arrow from
$x$ to $x+1$, then they exchange positions.  Otherwise, the
interaction is the usual exclusion.

Denote $X^k_t$ the position at time $t$ of the $\xi^k$ particle which
was at site $[c_k\ve^{-1}]$ at time $0$. 

The essential tools in [FF] (with $\mu_{\a,\be}$ as the initial
measure) include the joint realization of the
evolutions $(\eta'_t,Z_t)$ and $(\sigma^0_t,\xi^1_t,X^1_t)$, where
$\eta'_0$ is distributed as $\mu_{\a,\be}$ conditioned to 
have the site
$0$ occupied by a second class particle, and $Z_t$ describes the
position of this single second class particle, while 
$\sigma^k_t$, $\xi^k_t$ and $X^1_t$ were defined above
for $n=1$, $\rho_0=\a$, $\rho_1=\beta$. 

When $p=1$ the above coupling can be achieved by 
$$
\eta'_t(x) = \cases
\sigma^0_t(x)& \text{ if } x<X^1_t\\
\sigma^1_t(x)& \text{ if } x\ge X^1_t,\endcases
$$ 
and $Z_t=X^1_t$.  Together with the law of large numbers and the
central limit theorem for $X^1_t$, this coupling is the basic
ingredient in [FF]. Still restricting ourselves to the case $p=1$, a natural extension of this to our evolution
$\eta_t$ starting with $\mu^\ve$ is the following.  \vskip 3mm

\noindent{\bf Definition of $(Y^k_t)$.} 
Let us recall that we consider the case $p=1$, and $X^k_t$ are as just described above.
Let $Y^0_t\equiv-\infty$, $Y^{n+1}_t\equiv\infty$.
For
$k=1,\dots,n$, we define $Y^{k,i}_t$ and $t_i$ inductively in $i\geq0$ 
as follows.
Let $t_0=0$ and $Y^{k,0}_t=X^k_t$ for all $k=1,...,n$ and $t\geq0$. 
Having defined $t_\ell$ and 
$Y^{k,\ell}_t$ for all $k=1,...,n$ and $t\geq t_\ell$,
if $t_\ell=\infty$, we stop the inductive procedure; otherwise, let
$t_{\ell+1}= \inf \{t\geq t_{\ell}: Y^{i,\ell}_{t}=Y^{i+1,\ell}_{t}+1\}$ (with
the usual convention that $\inf\emptyset=\infty$). 
Denote by $i_\ell$ the index involved. If finite, 
$t_{\ell+1}$ is the time of the first crossing after time $t_\ell$
of two particles whose
positions at time $t_\ell$ are $Y^{\cdot,\ell}_{t_\ell}$. 
At time $t_{\ell+1}$,
a $\sigma^{i_\ell-1}|\sigma^{i_\ell+1}$ discrepancy appears at the position $Y^{i_\ell+1}_{t_{\ell+1}}$.
Since $p=1$, from $t_{\ell+1}$ on, $Y^{i_\ell,\ell}_{t}$ and $Y^{i_\ell+1,\ell}_{t}$ do not
uncross, so we can  ignore the
$\xi^{i_\ell}$ particles and consider only the evolution of $(\sigma^i:\,{i
\neq i_\ell})$. Accordingly, we define, for $t\geq t_{\ell+1}$, 
$Y^{i_\ell,\ell+1}_{t}=Y^{i_\ell+1,\ell+1}_{t}$ to be the position of this
$\sigma^{i_\ell-1}|\sigma^{i_\ell+1}$ discrepancy at time $t$. We also make
$Y^{i,\ell+1}_{t}=Y^{i,\ell}_{t}$ for $i\ne i_\ell,i_\ell+1$.
This concludes the inductive step. 

Notice that the inductive steps in the above
procedure correspond to instants of crossings of shocks, after each of which
we coalesce the crossing shocks into a single one, which is then made to
follow a new discrepancy, leaving the remaining shocks untouched.

We are ready to define  $(Y^k_t)$. For all $1\leq k\leq n$ and $t\geq0$, let
$Y^k_t=Y^{k,\ell}_{t}$, 
if
$t_\ell\leq t<t_{\ell+1}$. In case $t_n<\infty$, $t_{n+1}$ is defined
as $\infty$.

We further define
$$
\eta'_t(x)= 
\sum_{k=0}^n \sigma^k_t(x) \one\{Y^k_t\le x< Y^{k+1}_t\}.\teq(2.4)
$$

Notice that some of the indicator functions may vanish. This happens
when the corresponding $Y$ particles coalesce before $t$.  Since $p=1$,
$\eta'_t$ has the same law as the evolution starting from $\mu^\ve$ 
conditioned to have the sites $[c_k\ve^{-1}]$ occupied by second class
particles with respect to the other particles of $\eta'$ such that the 
second class particles  of lower labels have
priority over the ones with higher labels.

\medskip

\ppclaim \Theorem (t121). Writing $\tilde X^k_{\tau_\ve}=X^k_{\tau_\ve}- [x_\ve]$
and $\tilde Y^k_{\tau_\ve}=Y^k_{\tau_\ve}- [x_\ve]$, where
$\tau_\ve$ and $x_\ve$ are given by \equ(tre), we
have:
$$
\align
\lim_{\ve\to 0}
\ve^{1/2}( \tilde X^1_{\tau_\ve},\dots,\tilde X^n_{\tau_\ve})\ 
&\overset{\Cal D}\to{=}\ (\x_1,\dots,\x_n)\teq(p7) \cr
\lim_{\ve\to 0}\ve^{1/2}(\tilde Y^1_{\tau_\ve},\dots,\tilde Y^n_{\tau_\ve})\
&\overset{\Cal D}\to{=}\ (\y_1,\dots,\y_n)\teq(p8)
\endalign 
$$ 
where $(\x_k)$ are i.i.d. centered Gaussian with variances $D_k$ given
by \equ(p99) and $(\y_1, \dots, \y_n) =
\psi(\x_1,\dots,\x_n)$, with $\psi$ defined in Section 1.

Before proving Theorem \equ (t121), let us recall that from Theorem 1.1 of
[FF] we have
$$ 
\ve^{1/2}\tilde X^k_{\tau_\ve}\;\overset{\Cal D}\to{\underset{\ve
\to 0} \to{\longrightarrow}} \x_k,\teq(3.1)
$$ 
 with $(\tilde X^k_{\tau_\ve})$ and $(\x_k)$ as in the previous
statement. We discuss independence below.  The basic idea to show this
theorem is to look at the system ``just before'' $\tau_\ve$ (in
macroscopic time); more precisely: let $\a>0$ and
$\tau_\ve^\a=\tau_\ve-\a\ve^{-1/2}$. We shall make $\ve\to0$ and then
$\a\to+\infty$.

Presumably, if $\ve$ is small, but $\a$ is large enough, with
overwhelming probability the second class particles in the $n$
one-shock systems, $(X^1_{\tau_\ve^\a},\dots,X^n_{\tau_\ve^\a})$ have
not yet crossed each other and so
$Y^k_{\tau_\ve^\a}=X^k_{\tau_\ve^\a}$ for $k=1,\dots,n$. This is due to
the control on the asymptotic behavior of $X^k_t$, known from [FF] and
a reasoning as the one used to define $b_k(\un x)$ in \equ(p223). This
allows to compare the system with $n$ independent systems of only one
shock each, and we conclude that
$(Y^1_{\tau_\ve^\a},\dots,Y^n_{\tau_\ve^\a})$  are at
asymptotically independent Gaussian distances from their expected
values (on the $\ve^{-1/2}$ scale). (Notice that since $p=1$, then
each pair $X^i$ and $X^j$ with $i <j$ cross each other only once.)

Since, furthermore, $(Y^1_{\tau_\ve^\a},\dots,Y^n_{\tau_\ve^\a})$ are
at distance of order $\ve^{-1/2}$ apart from each other, in the time
interval from $\tau_\ve^\a$ to $\tau_\ve$ (of length $\a\ve^{-1/2}$)
the fluctuations are negligible (we are in the regime of the
hyperbolic scaling, thus of the law of large numbers for $Y^k_\cdot$)
and it is as if we have a completely deterministic situation.

Rigorous versions of the facts described on the previous two paragraphs 
lead straightforwardly to Theorem \equ(t121).
We state and prove them in the paragraphs below.

\bigskip

The first step is to describe the behavior of
$(X^1_{\tau_\ve^\a},\dots,X^n_{\tau_\ve^\a})$. Since each of the
coordinates $X^k$ correspond to the motion of a second class tagged
particle in a single shock situation, the essential is their
dependence on the initial condition. Indeed, if
$$ 
N_{t,\ve}^k(\sigma^{k-1},\sigma^k) \overset{\text{def}}\to{=}
\sum_{x=[c_k/\ve]+1}^{[(c_k+t(\rho_{k}-\rho_{k-1}))/\ve]} (1-
\sigma^k(x))-\sum_{x=[(c_k-t(\rho_{k}-\rho_{k-1}))/\ve]+1}^{[c_k/\ve]}\sigma^{k-1}(x)
$$
it follows from [FF] that 
$$ 
\ve t^{-1}E\lv X^k_{t/\ve}-c_k-
\frac{N^k_{t,\ve}(\sigma^{k-1},\sigma^k)}{\rho_k-\rho_{k-1}} \rv^2
\underset{t\to+\infty}\to{\longrightarrow} 0\teq(3.3a)
$$
  
\noindent {\bf Remark.} Indeed the result in [FF], for $X^k_t$, is for
uniform product measures off the initial location of the tagged
particle. A straightforward conditioning argument shows it to be true
for any finite perturbation of that measure. From Theorem 3.1 in
(Ferrari (1992)), but also from \equ(3.3a) it holds
$$
\lim_{t\to\infty} {\ve X^k_{\ve^{-1}t}} \,=\,
c_k+(1-\rho_{k-1}-\rho_k)t,\;\;\;\text{ in probability}. \teq(f31)
$$ 

\noindent{\bf Proof of \equ(p7).}  If $t\le t^*$ the random variables
$N^k_{t,\ve}$ are independent since by \equ(p10) they are functions of
the initial configuration at disjoint sets of sites. This and the
product nature of the initial distribution allows us to use the
central limit theorem for sums of Bernoulli random variables to show
\equ(p7). \square

To compute the limiting joint probabilities of \equ(p8) let $\Delta>0$
and consider the following intervals
$$ 
\eqalign{ J_i&:= \Delta[i-1/2,i+1/2)\cr J^\ve_i&:= \ve^{-1/2}J_i =
\left[\un\jmath_i,\ov\jmath_i\right)\;\;\text{ centered at }\;\;
\jmath_i := i \,\Delta\ve^{-1/2},\;\;\; i\in\Z.} 
$$

The following lemma takes care of the case when the $X^k_{\tau_\ve}$ have not still 
crossed each other.
\bigskip

\ppclaim\Lemma (l3.1). For integers $i(1)<\dots<i(n)$, 
$$
\lim_{\ve\to0}\,P(\tilde Y^k_{\tau_\ve}\in J^\ve_{i(k)}, k=1,\dots,n)
=\,\prod_k P\bigl(\x_k\in J_{i(k)}\bigr),
           \teq(3.4)
$$
where $(\x_k)_k$ are as in Theorem \equ(t121).

\flushpar {\bf Proof:} It follows from the definition of $Y^k_t$ and
\equ(p7). \square \bigskip

Now we study the case of unordered $X^k_{\tau_\ve}$. We consider first the
case $n=2$. Take arbitrary integers $i(1)>i(2)$. For $0<\delta<1/2$ let
$$ 
J^{\ve,\delta}_i:= \Bigl(i-\frac12+\delta,\,i+\frac12-\delta \Bigr)
\Delta\ve^{-1/2}\, =\, \left(\un\jmath_i+\delta\Delta\ve^{-1/2},
\,\ov\jmath_i-\delta\Delta\ve^{-1/2}\right). 
$$ 
The interval
$J^{\ve,\delta}_i$
is strictly contained in $J^\ve_i$.  Let
$$ 
I^{\ve,\delta}_{i(1),i(2)}\;=\;\Bigl(\psi_1(i(1),i(2))
-1-\delta,\,\psi_1(i(1),i(2))+1+\delta\Bigr) 
\,\, \Delta\ve^{-1/2}
$$
where
$$
\psi_1(i(1),i(2))\, =\, \psi_2(i(1),i(2))\, = \,
\displaystyle{i(1) {\rho_1-\rho_0\over \rho_2-\rho_0} + i(2)
{\rho_2-\rho_1\over \rho_2-\rho_0}}\teq(p29)
$$ 
This corresponds to the coordinates of the vector defined in
\equ(p28) because $i(1)>i(2)$.

\bigskip

\ppclaim\Lemma (l3.2). For integers $i(1)>i(2)$, and for all fixed $\Delta$, 
%there exists a positive constant $c$ such that
$$
\aligned &\lim_{\de\to0}\lim_{\ve\to0}\,
P\left(\tilde X^k_{\tau_\ve}\in J^{\ve,\delta}_{i(k)},\, k=1,2\right)\, 
=\,\lim_{\de\to0}\lim_{\ve\to0}\,
P\left(\tilde X^k_{\tau_\ve}\in J^{\ve,\delta}_{i(k)},\,
\tilde Y^k_{\tau_\ve}\in I^{\ve,\delta}_{i(1),i(2)},\,k=1,2\right).
\endaligned\teq(3.6)
$$

\bigskip

\flushpar
{\bf Proof:} 
For $I,J\subset \r$ we say 
$$
I<J \text\;\;{ if }\;\;x<y \;\;\text{ for
all\;\; $x\in I$\;\; and\;\; $y\in J$.}
$$
Considering, as before, $\tau_\ve^\a = \tau_\ve- \a \ve^{-1/2}$, let us  
define the intervals 
$$ 
\eqalign{ J^\ve_{i(k),\a}\,&=\,
(\un\jmath_{i(k)}^\a,\,\ov\jmath_{i(k)}^\a)\,=\,\theta_{-\a
\ve^{-1/2}(1-\rho_{k-1}-\rho_{k})\,+\,x_\ve} J^\ve_{i(k)},\cr
J^{\ve,\delta}_{i(k),\a}\,&=\, (\un\jmath_{i(k)}^\a+\delta
\Delta\ve^{-1/2},\,\ov\jmath_{i(k)}^\a-\delta
\Delta\ve^{-1/2})\,=\,\theta_{-\a
\ve^{-1/2}(1-\rho_{k-1}-\rho_{k})\,+\, x_\ve} J^{\ve,\delta}_{i(k)}}
$$ 
which we get at time $\tau_\ve^\a$ by translating
$\theta_{x_\ve}J^{\ve}_{i(k)}$ and $\theta_{x_\ve}J^{\ve,\delta}_{i(k)}$ 
at velocity $-(1-\rho_{k-1}-\rho_{k})$ backwards from
time $\tau_\ve$ to time $\tau^\a_\ve$. Let $\a$ big enough such that
$$
J^{\ve}_{i(1),\a}<J^{\ve}_{i(2),\a}.
$$ 
 From the one-shock law of large numbers for $X_t$ and the
definition of $Y^k_t$ (which coincide with the $X^k_t$ for $t < t_1$
in that defintion), we have:
$$
\align
\lim_{\ve\to 0}P(Y^k_{\tau_\ve^\a} \in J^\ve_{i(k),\a}, k=1,2) 
\,&=\,\lim_{\ve\to 0}P(X^k_{\tau_\ve^\a} \in J^\ve_{i(k),\a}, k=1,2)\cr
\,&=\, \lim_{\ve\to 0}P(\tilde X^k_{\tau_\ve}\in J^\ve_{i(k)}, k=1,2)\cr
\endalign
$$ 
This gives the localization of the $Y$ particles at time
$\tau^\a_\ve$.  The idea is that those particles will follow their
respective characteristics and thus meet during the time interval
$(\tau_\ve^\a,\tau_\ve)$, where they coalesce and follow the new
characteristic, thus ending in $x_\ve \,+\,I^{\ve,\delta}_{i(1),i(2)}$
at time $\tau_\ve$. Since we are in the same scale for time and space,
the result will follow from the law of large numbers.

To make this rigorous, at time $\tau_\ve^\a$ put $\sigma^0|\sigma^1$
discrepancies (second class particles) at 
$[\un\jmath_{i(1)}^\a+\delta \Delta\ve^{-1/2}]+1$ and 
$[\ov\jmath_{i(1)}^\a-\delta\Delta\ve^{-1/2}]$ and $\sigma^1|\sigma^2$
discrepancies (third class particles) at
$[\un\jmath_{i(2)}^\a+\delta\Delta\ve^{-1/2}]+1$ and
$[\ov\jmath_{i(2)}^\a-\delta\Delta\ve^{-1/2}]$.  Label their positions $\un
Y{}^1_t$, $\ov Y{}^1_t$, $\un Y{}^2_t$ and $\ov Y{}^2_t$,
respectively, for $t\geq\tau_\ve^\a$.

Let $\tau'_\ve=\tau''_\ve+\delta\ve^{-1/2}$, where  
$\tau''_\ve$ is defined by
$E\un Y^1_{\tau''_\ve}= E \ov Y^2_{\tau''_\ve}$.
At time $\tau_\ve'$, which will be smaller than $\tau_\ve$ for $\ve$ 
sufficiently small, put
$\sigma^0|\sigma^2$ discrepancies at 
$$
[\ov\jmath_{i(1)}^\a+(1-\rho_{0}-\rho_1)(\tau_\ve'-\tau_\ve^\a)]
\quad\text{and}\quad
[\un\jmath_{i(2)}^\a+(1-\rho_{1}-\rho_2)(\tau_\ve'-\tau_\ve^\a)].
$$ 
Label their positions $\un Y_t$ and $\ov Y_t$, 
respectively, for $t\geq\tau_\ve'$.
%This implies that for $t \geq \tau'_\ve$,
%$E(\ov Y_t)- E(\ov Y^k_t)$  and  $E(\un Y^k_t)- E(\un Y_t)$ differ from
%$\delta\ve^{-1/2}$ by at most $1$. 

Since $p=1$, we easily check the facts that for $t \geq \tau^\a_\ve$,
$$
\left\{Y^k_{\tau_\ve^\a}\in
J^{\ve,\de}_{i(k),\a},\,k=1,2\right\} 
=\left\{X^k_{\tau_\ve^\a}\in
J^{\ve,\de}_{i(k),\a},\,k=1,2\right\} \subseteq 
\left\{\un Y^k_t\leq X^k_t\leq\ov Y^k_t \right\}
$$
which shows in particular that for all $t \geq \tau'_\ve$: 
$$
\left\{Y^k_{\tau_\ve^\a}\in J^{\ve,\de}_{i(k),\a},\,k=1,2\right\} 
\subseteq \left\{ X^2_t \le X^1_t \right \}
\cup \left\{\un Y^1_{\tau'_\ve}\leq \ov Y^2_{\tau'_\ve}\right\} 
\teq(p44)
$$
and 
$$ 
\left\{ X^2_t\leq X^1_t\right\} \subseteq \left\{X^2_t\leq
Y^k_t\leq X^1_t,\,k=1,2 \right\}.\; \teq(p45)
$$

Now, a simple geometric argument relying on the laws of large numbers
for $\un Y{}^k_t, \ov Y{}^k_t, X^k_t$ (cf. Theorem 3.1 in [FF] and
Remark following eq.~\equ(3.3a)), \equ(p44) and \equ(p45) proves that
$$ 
\lim_{\ve\to0}\,P\left(Y^k_{\tau_\ve^\a}\in
J^{\ve,\de}_{i(k),\a},\,k=1,2\right)
=\lim_{\ve\to0}\,P\left(X^k_{\tau_\ve^\a}\in
J^{\ve,\de}_{i(k),\a}, \,\,\un
Y_{\tau_\ve'}\leq Y^k_{\tau_\ve'}\leq\ov Y_{\tau_\ve'},\,k=1,2\right)
\teq(3.7)
$$
plus $O(\delta)$, as soon as $\de>0$ is small enough. $O(\delta)$
comes from the probability of the event on the right of the union
on the right hand side  of  eq.~\equ(p44) and the law of large numbers 
for $\un Y^1_t$ and $\ov Y^2_t$.

The result now follows from an application of the laws of large numbers for 
$\un Y_{t}$, $\ov Y_{t}$ and the fact 
$\un Y_{t}\leq Y_{t}\leq\ov Y_{t}$ for $t\geq\tau_\ve'$. 
(See figure 3.1.) \square

\bigskip
\centerline{\hbox{\tt (Figure 3.1 enters here)}}
\bigskip

The same argument with a more complicated notation shows the above
extends to the case of general $n$, summarized below.  \bigskip

\ppclaim\Lemma (lpsi). For distinct integers $i(1),\dots, i(n)$, 
and for all fixed $\Delta$, 
$$
\aligned &\lim_{\de\to0}\lim_{\ve\to0}\,
P\left(\tilde X^k_{\tau_\ve}\in J^{\ve,\delta}_{i(k)},\, k=1,\dots,n\right)\cr 
&\qquad=\lim_{\de\to0}\lim_{\ve\to0}\,
P\left(\tilde X^k_{\tau_\ve}\in J^{\ve,\delta}_{i(k)},\,
       \tilde Y^k_{\tau_\ve}\in I^{\ve,\delta}_{j(k)},\,k=1,\dots,n\right),
\endaligned\teq(3.6p)
$$
where $j(k) = \psi_k(i(1),\dots,i(n))$.

\bigskip

\noindent{\bf  Proof of \equ(p8).} Follows from Lemma
\equ(lpsi) and \equ(p7). \square

\bigskip
\flushpar {\bf \S 3.\quad Proof of Theorems \equ(t1.1) and \equ(t1.2).}
\nopagebreak
\numsec=3\numfor=1

\bigskip

\flushpar{\bf Proof of Theorem \equ(t1.2).} 
Let us start by noticing that Theorem \equ(t1.2) was stated in terms of the
evolution $\eta_t$, which is not exactly the $\eta'_t$ we are considering, for which initially the sites $[c_k\ve^{-1}],\,k=1,\dots,n$, are occupied by second class particles. On the other side, it is immediate to see that  
coupling $\eta_t$ and 
$\eta'_t$ in the usual way, one concludes that the eventual discrepancies, 
at most $n$, behave as second class particles. They might even annihilate 
one another, but in any case diffuse as $\ve^{-1/2}$. This shows that
Theorem \equ(t1.2) will be proven (i.e. for $\eta_t$) once we check the
analogous for our perturbed process $\eta'_t$. For this,  
it is sufficient, according to Daley and Vere-Jones
(1988), Proposition 9.1.VII, to show that for every bounded continuous
function $\Phi$ with compact support $\int \Phi\,d\La_\ve$ converges
weakly to $\int\Phi\,d\La$ as $\ve\to0$. The former integral is
$$
\ve^{1/2}\sum_{x\in\Z}f(\theta_{x+[x_\ve]}\eta'_{\tau_\ve})\,\Phi(\ve^{1/2}x).
\teq(5.3) 
$$

Let $M \geq 0$ be an integer such that the cylinder function $f$ depends only on the coordinates in $\{-M,\dots,M\}$. 
Then, we use \equ(2.4) to condition on the (standardized) locations
of the shocks at time $\tau_\ve$, and decompose \equ(5.3) as
$$ 
\sum_{k=0}^{n}
\ve^{1/2}\sum_{x\in\Z}f(\theta_{x+[x_\ve]}\si^k_{\tau_\ve})\,\Phi(\ve^{1/2}x)
\,\one\{\hat Y_{\ve}^{k}\le\ve^{1/2}(x-M) <\ve^{1/2}(x+M) <\hat
Y^{k+1}_{\ve}\}\,+\, Z_\ve \teq(5.4)
$$
where
$\hat Y^k_{\ve}=\ve^{1/2}(Y^k_{\tau_\ve}-[x_\ve])$ and the random variable
$Z_\ve$ satisfies 
$$
E(|Z_\ve|)\leq ||f||_\infty\,\sum_{k=0}^{n}\ve^{1/2}\sum_{x\in\Z}\,\Phi(\ve^{1/2}x)P\left(\hat Y^k_{\ve} \in (\ve^{1/2}(x-M),\ve^{1/2}(x+M)\right). 
$$

By Theorem \equ(t121), $(\hat Y^1_{\ve},\dots,\hat Y^n_{\ve})$ converges weakly
to $\psi(\x_1,\dots,\x_n)$ as $\ve\to0$, and this implies at once that $\lim_{\ve \to 0}E(|Z_\ve|)=0$.

The result follows from this. To see it, notice that,
for $\ve>0$, the expression in \equ(5.4) is a
function of $(\si^k_{\tau_\ve})$ and $(\hat Y^k_{\ve})$, 
which we denote
$F_\ve(\si^1_{\tau_\ve},\dots,\sigma^n_{\tau_\ve}, \hat
Y^1_{\ve},\dots,\hat Y^n_{\ve})$.  By the product structure of
$\si^k_{\tau_\ve}$ and the law of large numbers, for all
$(y_1,\dots,y_n)$
$$
F_\ve(\si^1_{\tau_\ve},\dots,\sigma^n_{\tau_\ve},
y_1,\dots,y_n) \,\underset{\ve\to0}\to{\longrightarrow}\, 
F(y_1,\dots,y_n):=\int
\Phi\,d\La_{(y_1,\dots,y_n)}
$$ 
almost surely, where $\La_{(y_1,\dots,y_n)}$ is defined as in
\equ(5.2). Both $F$ and $\psi$ are continuous in $(y_1,\dots,y_n)$. 
Also, $\{F_\ve,\, \ve>0\}$ is equicontinuous as a family of
functions of $(y_1,\dots,y_n)$.  This almost sure convergence and
continuity properties of $F_\ve$, $F$ and $\psi$ plus the weak
convergence of $(\hat Y^k_{\ve})$ yield the weak convergence claimed
in the statement of the theorem by standard arguments.  We leave the
details to the reader. \square

\bigskip

\flushpar {\bf Proof of Theorem \equ(t1.1).} It suffices to check that 
$$
\lim_{\ve \to 0}E f(\theta_{[r^*\ve^{-
1}+a\ve^{-1/2}]}\eta_{t^*\ve^{-1}})\ \overset{\om^*}\to{\underset{\ve\to0}\to{\longrightarrow}}\
\sum_{k=0}^{n}    \nu_{\rho_k}f\,P(\y_{k} \le r < \y_{k+1}) 
$$ 
where $\y_k$ are as in the statement of Theorem \equ(t1.1), $\eta_t$ represents the process starting with $\mu^\ve$ and $f$ is a
cylinder function, increasing (for the usual coordinatewise order).
Now, to recover the above expression out of Theorem
\equ(t1.2) we may take the expectation $E \int \Phi d\Lambda_\ve$
where the test function $\Phi$ defined by $\Phi(w)=\frac {1}{u}
\one\{a\le w<a+u\}$, with $u>0$ will give an upper bound for $E
f(\theta_{[x_\ve + a\ve^{-1/2}]} \eta_{\tau_\ve})$. Similarly by using $\Phi$
defined by $\Phi(w)=\frac {1}{u} \one\{a-u\le w<a\}$, with $u>0$, we
get a lower bound. (The fact that these particular test functions have two
points of discontinuity is not a problem, due to the continuity of limit
measure $\Lambda$.)
Here we are using the attractiveness of the system
and the monotonicity of the initial profile to get $\mu^\ve S(t)f \le
\theta_1\mu^\ve S(t)f$, for increasing continuous $f$.  Letting $u$ tend to zero
both terms will converge to the desired expression. \square

\medskip

\
\flushpar {\Remark(3.3)} A similar analysis can be employed to
determine that the (microscopic) measure seen from
($[x_\ve+a\ve^{-1/2}],\tau_\ve+s\ve^{-1/2}$), with $a,s\in\r$ fixed,
converges as $\ve\to0$ to a mixture of $\nu_{\ro_k}$. For that, as in
the definition of $\psi$, let 
$$
\psi^s_k(\un x)=
\cases  x_k + s (1-\rho_{k-1}-\rho_k) &\text{if }s\le -t(\un x)\cr
b_k(\un x, t(\un x)+s)&\text{if }s> -t(\un x)\cr
\endcases
$$ 
This is the function that enters in the corresponding statement.

\medskip

\flushpar {\Remark(3.5)} In this paper we have treated in detail only
the totally asymmetric case, where $p=1$. It is not hard, but
technically more cumbersome, to extend the analysis to $1/2<p<1$. One
can define $\eta'_t$ in the same way. The extra difficulty comes from
the fact that now the shocks can {\it un}cross each other. However,
since the gaussian fluctuations and the law of large numbers remain
valid, essentially the same analysis applies, with the appropriate
change on parameters. Details are left to the reader.

%%%%%%%%%%%%%%%%%%%%%%%%%%%%%%%%%%%%%%%%%%%%%%%%%%%%%%%%%%%%%%%%%%%%%%%%%

\bigskip
\noindent{\bf Acknowledgments.} 
This paper was partially supported by CNPq, FAPESP and N\'ucleo de Excel\^encia
``Fen\^omenos Cr\'\i ticos em Probabilidade e Processos Estoc\'asticos''.
We also acknowledge support from ProInter-USP. The final version was
written while the first author was invited professor at Laboratoire
de Probabilit\'es de l'Universit\'e Pierre et Marie Curie, Paris
VI. He thanks warm hospitality.

\bigskip

\noindent{\bf References.}

\item {--} E.\ D.\ Andjel, M.\ Bramson, T.\ M.\ Liggett (1988). Shocks in
the asymmetric simple exclusion process. {\sl Probab.\ Theor.\ Rel.\ Fields
\bf 78 } 231-247.

\item{--} E.\ D.\ Andjel, M.\ E.\ Vares (1987). Hydrodynamic equations for
attractive particle systems on $\Z$. {\sl J.\ Stat.\ Phys.\ \bf 47}
265-288.

\item {--} D.\ J.\ Daley, D.\ Vere-Jones (1988).
{\sl An Introduction to the Theory of Point Processes}.
Springer, Ber\-lin.

\item {--} A.\ De Masi, C.\ Kipnis, E.\ Presutti, E.\ Saada (1988).
Microscopic structure at the shock in the asymmetric simple exclusion. {\sl
Stochastics \bf 27}, 151-165.

\item {--} P.\ A.\ Ferrari (1992) Shock fluctuations in asymmetric simple
exclusion. {\sl Probab.\ Theor.\ Related Fields. \bf 91}, 81--101.

\item{--} P.\ A.\ Ferrari, L.\ R.\ G.\ Fontes (1994) Shock fluctuations for
the asymmetric simple exclusion process {\sl Probab. Theory
Related Fields \bf 99}, 305-319.

\item{--} P.\ A.\ Ferrari, C.\ Kipnis, E.\ Saada (1991) Microscopic
                  structure of travelling waves in the 
                  asymmetric simple exclusion process {\sl
                  Ann. Probab. \bf 19} 1:226-244.

\item {--} C. Landim (1992). Conservation of local equilibrium for
attractive particle systems on $Z^d$. {\sl Ann. Probab.}

\item {--} T.\ M.\ Liggett (1976). Coupling the simple exclusion process.
{\sl Ann.\ Probab.\ \bf 4} 339-356.

\item {--} T.\ M.\ Liggett (1985). {\sl Interacting Particle Systems.}
Springer, Ber\-lin.

\item {--} H.\ Rezakhanlou (1990) Hydrodynamic limit for attractive particle
systems on $Z^d$. {\sl Comm.\ Math.\ Phys. \bf 140} 417-448. 

\item {--} F.\ Spitzer (1970). Interaction of Markov processes. {\sl Adv.
Math., \bf 5} 246-290.

\item {--} D. Wick (1985). A dynamical phase transition in an infinite
particle system. {\sl J.\ Stat.\ Phys.\ \bf 38} 1015-1025.

\enddocument